\documentclass{amsart}

\usepackage{amssymb,amsfonts}

\usepackage[dvips]{graphicx}

\theoremstyle{plain}

\newtheorem{thm}{Theorem}[section]

\newtheorem{cor}[thm]{Corollary}

\newtheorem{example}[thm]{Example}

\theoremstyle{definition}

\theoremstyle{remark}

\newcommand{\Sth}{\textbf{S}^3}

\newcommand{\Hth}{\textbf{H}^3}

\begin{document}


\title[Totally Geodesic Surfaces in Knot and Link Complements I]
{Totally Geodesic Seifert Surfaces in Hyperbolic Knot and Link Complements I}

\date{\today}
\author[Colin Adams]{Colin Adams}
\author[Eric Schoenfeld]{Eric Schoenfeld}

\address{Colin Adams, Department of Mathematics and Statistics, Williams
College, Williamstown, MA 01267}
\email{Colin.C.Adams@williams.edu}

\address{Eric Schoenfeld, Department of Mathematics, Stanford University,
450 Serra Mall, Bldg. 380, Stanford, CA 94305-2125}
\email{erics@math.stanford.edu}

\begin{abstract}
The first examples of totally geodesic Seifert surfaces are constructed for
hyperbolic knots and links, including both free and totally knotted
surfaces. Then it is proved that two bridge knot complements cannot 
contain totally geodesic orientable surfaces.
\end{abstract}

\maketitle

\section{Introduction}\label{S:intro}
Define a knot or link in $\Sth$ to be hyperbolic if its complement is a
hyperbolic $3$-manifold. This implies that there is a covering map $p$
from $\Hth$ to $\Sth-K$ such that the covering translations are isometries
of $\Hth$. We say that an embedded or immersed surface $S$ in $\Sth-K$ is
totally geodesic if it is isotopic to a surface that lifts to a set of
geodesic planes in $\Hth$.

Embedded totally geodesic surfaces are relatively rare in knot and
link complements. It was conjectured in ~\cite{MR92} that there are no
closed embedded totally geodesic surfaces in hyperbolic knot complements.
This has been verified for hyperbolic knots that are alternating knots,
tunnel number one knots, and 2-generator knots ~\cite{MR92}, almost
alternating knots ~\cite{Aetal92}, toroidally alternating knots
~\cite{cA94}, Montesinos knots ~\cite{uO84}, 3-bridge knots and double
torus knots ~\cite{IO200}) and knots of braid index three ~\cite{LP85} and
four ~\cite{hM02}. Note that the conjecture does not hold for links.  The
first explicit counterexample was given in ~\cite{MR92}. In  ~\cite{Le04},
it is demonstrated that for each $g > 2$, there exists a 2-component
hyperbolic link that contains a closed embedded totally geodesic genus $g$
surface.

In this paper, we investigate surfaces that are not closed in knot and
link complements. In particular, we are interested in totally geodesic
Seifert
surfaces, where a Seifert surface for a knot or link $L$ is a compact
orientable
surface whose boundary is $L$. Note that although a totally geodesic
Seifert surface need not necessarily be minimal genus, it is
incompressible.

Work of Thurston implies that a properly embedded surface $S$ with
boundary in a hyperbolic knot exterior $M$ can have one of three
possible behaviors:

\begin{enumerate}

\item There is a nontrivial simple closed curve $c$ in $S$ that is
not boundary parallel in $S$ but that corresponds to a parabolic isometry.
This implies that $c$ is isotopic into a neighborhood of the missing knot.
In this case, we say that $c$ is an accidental parabolic curve and $S$ is
an accidental surface. See ~\cite{IO200} for more on these surfaces.

\item $S$ is a quasi-Fuchsian surface, with limit set a union of
quasi-circles. This includes the case of $S$ being totally geodesic (also
called Fuchsian), in which case the limit set is a union of geometric
circles.

\item $K$ is a fibered knot and in some finite cover, $S$ lifts to a
fiber. In this case, the limit set is the entire sphere.

\end{enumerate}

In ~\cite{sF98}, it is proved that if $K$ is not fibered, then any minimal
genus Seifert surface cannot be accidental and hence must be
quasi-Fuchsian. See also  ~\cite{CL01}.A quasi-Fuchsian surface $S$ 
can admit an essential annuli
$A$
intersecting $S$ in its boundary. However, in this case, $S$ cannot be
totally geodesic. In essence, the existence of the
annulus forces pairs of quasi-circles in the limit set of the cover to
share two points, which a pair of circles that do not cross cannot do. So,
for a surface to be totally geodesic, it must not admit any essential
annuli. Nor can it admit an essential disk with boundary consisting of four
arcs, two of which lie in the surface and two of which lie in the cusp.
Such a disk doubles to an essential annulus in the manifold obtained by
doubling along the totally geodesic surface.  However this contradicts the
fact the double is a hyperbolic manifold.

A Seifert surface is said to be {\em free} if the complement of a
neighborhood of the Seifert surface is a
handlebody. Note that there are examples of knots such that a minimal genus
Seifert surface is free and examples such that it is not free. See
~\cite{yM87}.

In Section \ref{S:tgss}, we explain how lifts of rigid hyperbolic
2-orbifolds can be utilized to generate totally geodesic Seifert
surfaces.
In Section \ref{S:examples}, we utilize these results to generate the
first examples of hyperbolic knots with totally geodesic Seifert surfaces.
We include both examples of free totally geodesic Seifert surfaces and
non-free totally geodesic Seifert surfaces. In Section \ref{S:twobridge},
we prove that many hyperbolic knots cannot
possess a totally geodesic Seifert surface, including all two-bridge
knots.

Note that the surfaces given here are the only known examples of embedded
totally geodesic surfaces of any kind in knot complements in $\Sth$. One
can ask whether there are embedded totally geodesic
surfaces in knot complements other than Seifert surfaces.

\medskip

Acknowledgements:
Thanks to Satyan Devadoss, Cameron Gordon and Alan Reid
for helpful conversations.

\section{Producing Totally Geodesic Seifert Surfaces}\label{S:tgss}

We call a hyperbolic orbifold {\em rigid} if there exists only one
possible hyperbolic metric that generates the orbifold. In ~\cite{wT78},
it is
demonstrated that a hyperbolic 2-orbifold $Q$ with underlying space $X_Q$,
$j$ elliptic points and $k$ corner reflectors is rigid if and only if
$3\chi(X_Q) = 2j + k$. A 2-orbifold is hyperbolic if and only if its Euler
characteristic is negative where the Euler characteristic is given by:

$$\chi(Q) =\chi(X_Q) - \sum(1-\frac{1}{m_i}) -
\frac{1}{2}\sum(1-\frac{1}{n_i})$$
where ${m_1, m_2, \dots, m_j}$ and ${n_1, n_2, \dots, n_k}$ are the orders
of the elliptic points and corner reflectors respectively.

Hence a complete list of rigid hyperbolic 2-orbifolds consists of

\begin{itemize}

\item $S^2(m_1, m_2, m_3;)$ for $1/m_1 + 1/m_2 + 1/m_3 < 1$,

\item $D^2(;n_1,n_2,n_3)$ for  $1/n_1 + 1/n_2 + 1/n_3 < 1$

\item $D^2(m_1; n_1)$ for  $2/m_1 + 1/n_1  < 1$

\end{itemize}

Note that each $n_i$ and $m_i$ can be either a positive
integer or $\infty$.

The fundamental group $\pi_1(O)$ of an orbifold $O$
is the group of deck transformations of the universal cover $\tilde{O}$.
A $2$-orbifold $Q$ embedded in a $3$-orbifold $O$ is {\em essential} if
$\pi_1(Q)$ injects into $\pi_1(O)$. A hyperbolic 2-orbifold $Q$ embedded
in a hyperbolic 3-orbifold $O$ is {\em totally geodesic} if it lifts to the
disjoint union of geodesic planes in the universal cover $\Hth$. Let $p$
denote the map projecting $\Hth$ onto $O$. We state a theorem that follows
  from the theorem on p. 217   ~\cite{bM88}, but provide a complete proof here
  for convenience.

\begin{thm}\label{main}
Let $Q$ be an essential  rigid hyperbolic $2$-orbifold embedded in a
hyperbolic $3$-orbifold $O$.  Then $Q$ is isotopic relative the singular
set to a totally geodesic $2$-orbifold $Q'$.
\end{thm}

\begin{proof}
Every isometry $\gamma$ of hyperbolic 3-space can be represented as
the product of two order two elliptic isometries, each a rotation about a
geodesic axis. If $\gamma$ is hyperbolic, then the two axes, considered at
subsets of $\Hth \cup \partial \Hth$,  do not intersect. If $\gamma$ is
parabolic then they share an endpoint. If $\gamma$ is elliptic, then they
intersect in $\Hth$ and the axis of rotation for the elliptic isometry
passes through the point of intersection and is perpendicular to the plane
containing both of the axes.

The fundamental group $\pi_1(O)$ is realized as a discrete group
$\Gamma$ of isometries of $\Hth$  with subgroup $\Gamma'$ corresponding to
$\pi_1(Q)$. There exists a pair of isometries  $g_1$ and $g_2$ in
$\Gamma'$, the product of which is a third isometry $g_3$ in $\Gamma'$. In
the cases of $S^2(m_1, m_2, m_3;)$ and $D^2(;n_1,n_2,n_3)$, we can choose
these isometries corresponding to the three singularities and the three
corner reflectors respectively.  In the case of $D^2(m_1; n_1)$, we can
choose isometries corresponding to the corner reflector and to two lifts
of the isometry corresponding to the singularity, one of which is conjgate
to the other by an order two elliptic isometry corresponding to the
boundary of the orbifold.

Each $g_i$ can be realized as the product of two isometries, each of which
is an order two elliptic isometry about an axis. Choose an axis for an
order
two elliptic isometry $f_2$ such that if $g_1$ is elliptic, the axis is
perpendicular
to the axis of $g_1$. If $g_1$ is parabolic, choose the axis of $f_2$ so
that it  ends at the
fixed point of $g_1$. Also make the axis simultaneously satisfy the same
requirements for $g_2$. Then we can find order two elliptic isometries
$f_1$ and $f_3$ such that $g_1 = f_1 f_2$ and $g_2 = f_2 f_3$. Note then
that $g_3 =g_1 g_2 = f_1 f_2 f_2 f_3 = f_1 f_3$. Since
$g_3$ is elliptic or parabolic, the axes for $f_1$ and $f_3$ must either
intersect in $\Hth$ or in $\partial \Hth$. In either case, the three
axes for $f_1$, $f_2$ and $f_3$ all lie in a geodesic plane $P'$.  The
entire group $\Gamma'$ is generated by $f_1$, $f_2$ and $f_3$ so $\Gamma'$
preserves $P'$. The orders of $g_1$, $g_2$  and $g_3$
determine the angles of the triangle defined by the
axes of  $f_1$, $f_2$ and $f_3$, and hence the resulting 2-dimensional
hyperbolic orbifold $Q'$ that is the projection of $P'$ to $O$ is rigid.

Exactly as in the proof of Theorem 3.1 of ~\cite{cA85}, we can prove that
the limit set of $P'$ is in fact also the limit set of a topological
plane $P$ covering $Q$. Since $Q$ is embedded, the limit sets of the
topological planes covering $Q$ can only intersect tangentially. Hence the
same holds for all the copies of $P'$ under the group action, implying
that the geodesic planes covering $Q'$ do not intersect in their
interiors. In particular, this means that $Q'$ is embedded and
$p^{-1}(Q')$ is a disjoint union of geodesic planes, and $Q'$ is a totally
geodesic orbifold.

Suppose $\mu$ is an element of $\Gamma$ such that
$\mu(P') = P'$. Then $\mu$ preserves the limit circle of $P'$. But then
$\mu$ preserves the limit set of $P$, implying $\mu \in \Gamma'$. Hence,
$Q$ and $Q'$ have the same fundamental group and are orbifold
homeomorphic.

Note that in the case the underlying space of $Q'$ is a
disk, $Q$ and $Q'$ must share the boundary arcs, which are order two
singularity axes in $O$.

The orbifold $Q$ can be isotoped so that it does not intersect $Q'$ in its
interior as follows. First, put $Q$ in general position relative to $Q'$.
Utilizing
the product neighborhoods of the singularities, cusps and boundaries of
$Q'$,
we can isotope $Q$ so that there are no intersections curves near
singularities,
cusps, or boundaries of the two 2-orbifolds. Hence, all the intersection
curves are simple closed curves. A trivial intersection curve on one of
the 2-orbifolds must be trivial on the other by the fact $Q$ is
essential and $Q'$ is totally geodesic. However, such a curve can then be
eliminated by isotoping $Q$.

A curve that is nontrivial on both $Q$ and $Q'$, and innermost such
on $Q$ must bound a disk with one singularity on each of $Q$ and $Q'$,
denoted $D$ and $D'$ respectively. Since this implies that a loop around
the singularity on $Q$ is homotopic to the loop around the singularity
on $Q'$, the singularities must either both correspond to parabolic
isometries or to elliptic isometries of the same order.  Then $D \cup D'$
forms a sphere with two singularities of the same order.  By
irreducibility, $D \cup D'$ bounds a ball containing an unknotted
singular arc or cusp within it. Hence, $D$ can be isotoped through $D'$
to eliminate the intersection.

Now, $Q$ and $Q'$ do not intersect. Hence their lifts, which come paired
with shared boundary circles, do not intersect. Thus each pair of planes
sharing a boundary circle bound a ball in $\Hth$ that is intersected by
none of its translates in the interior. Thus, we can isotope the
topological plane that is a lift of $Q$ to the corresponding geodesic
plane that is a lift of $Q'$. In fact, we can do so equivariantly with
respect to $\Gamma$, and none of the disjoint balls
bounded by the planes will intersect in the process. Hence $Q$ is isotopic
to $Q'$, as we wished to demonstrate.

\end{proof}

\begin{cor}\label{C:totallygeod}
Let $M$ be a closed orientable 3-manifold, covering a 3-orbifold $O$ such
that $M $ contains a link $L$ covering a collection $J$ of simple closed
curves and arcs beginning and ending on the order two singular set of $O$,
and
such that $O-J$ is a hyperbolic 3-orbifold $O'$. Then if $O'$ contains a
rigid
hyperbolic 2-orbifold $Q$ for which each simple closed component of $J$ is
a single puncture (singularity with $m_i=\infty$) and each arc corresponds
to a corner reflector with $n_j=\infty$, then $Q$ is isotopic in $O'$to a
2-orbifold that lifts to a totally geodesic Seifert surface for the link
$L$ in $M$.
\end{cor}

\begin{proof}
This follows immediately from Theorem \ref{main}.
\end{proof}

\section{Examples}\label{S:examples}

In this section, we utilize Corollary
~\ref{C:totallygeod} to give a variety of examples of totally geodesic
Seifert
surfaces in knot and link complements in $\Sth$.

\medskip

\begin{example}Balanced Pretzel Links
\end{example}

\begin{figure}[htbp]
\begin{center}
\includegraphics*[height=2in]{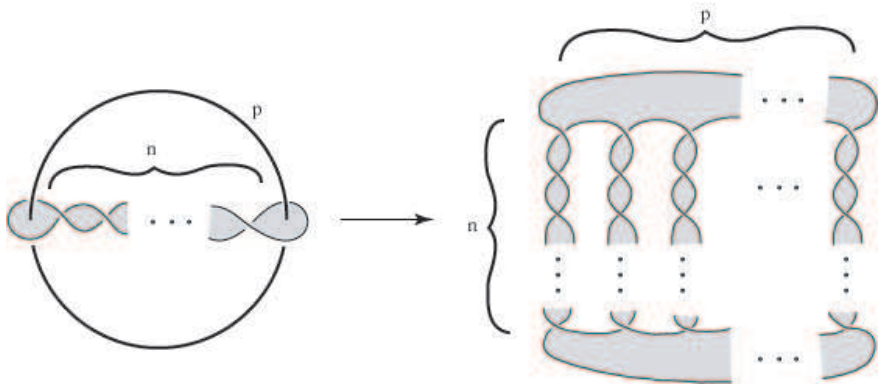}
\end{center}
\caption{\label{fig:balancedpretzel}The Seifert surfaces for this family
of balanced pretzel links are all totally geodesic.}
\end{figure}

\begin{figure}[htbp]
\begin{center}
\includegraphics*[height=4 in]{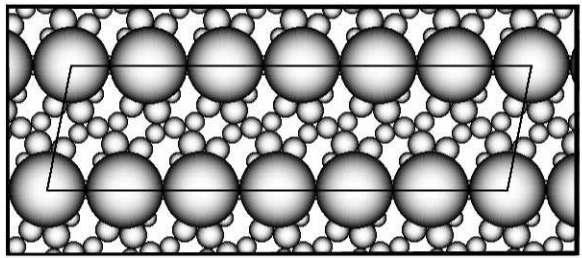}
\end{center}
\caption{\label{fig:balancedpretzellift}The straight line of horoballs
corresponds to the totally geodesic surface.}
\end{figure}

Let $M$ be $\Sth$ and $O$ the orbifold obtained by quotienting out by a
$2\pi/p$ rotation about a great circle $C$ in $\Sth$, where $p \ge 2$. Let
$D$ be a disk with two singularities of order $p$ and $p$, twisted with
$n$ crossings as in Figure \ref{fig:balancedpretzel}. Let  $J$ be the
boundary of the disk $D$. Then $J$ lifts to an $(n, n, n,\ldots,n)$
pretzel knot $K$, where $n$ appears $p$ times. We call $K$ a {\em balanced
pretzel knot}. A picture of the horoball diagram for the $(3,3,3)$ pretzel
knot appears in Figure \ref{fig:balancedpretzellift} . This picture was
generated by Jeffrey Weeks' SnapPea computer program, (cf. ~\cite{jW03}).
The
sequentially tangent horoballs in a line
correspond to a vertical geodesic plane coming up out of the xy-plane that
covers the totally geodesic surface.

\begin{example}Polyhedral links
\end{example}

\begin{figure}[htbp]
\begin{center}
\includegraphics*[height=2 in]{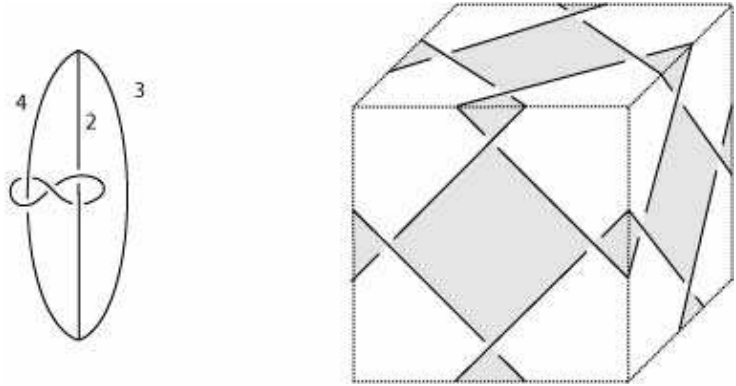}
\end{center}
\caption{\label{fig:cubelink}A totally geodesic Seifert surface for this
cube link.}
\end{figure}

A link that projects to the surface of the cube, as in Figure
\ref{fig:cubelink}, corresponds to the loop in the orbifold shown, which
possesses
the shaded rigid 2-orbifold that appears in its complement.
Therefore, the link bounds the totally geodesic Seifert surface shown.

\begin{example}Knots with Totally Knotted Seifert Surfaces
\end{example}

\begin{figure}[htbp]
\begin{center}
\includegraphics*[height=2 in]{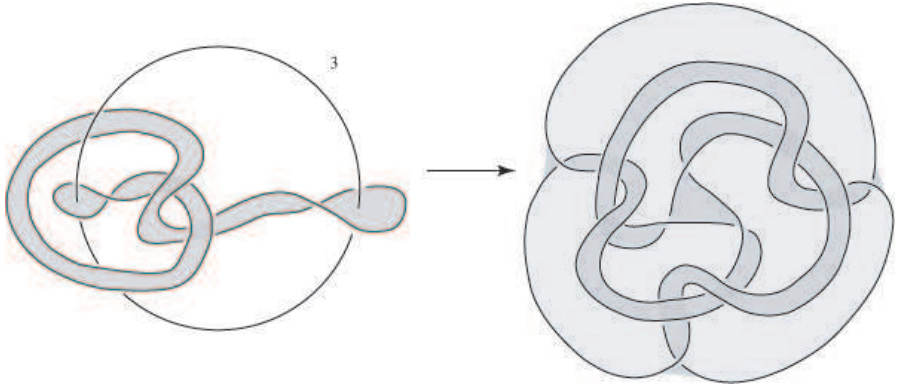}
\end{center}
\caption{\label{fig:totallyknotted}A totally knotted Seifert surface which
is totally geodesic.}
\end{figure}

Instead of a disk with two cone points that is twisted relative the
singular arc as we saw in Figure \ref{fig:balancedpretzel}, we will choose
a disk that is knotted relative the singular curve as in Figure
\ref{fig:totallyknotted}.
When $p = 3$, this lifts to the knot in Figure \ref{fig:totallyknotted}.
The disk lifts to a totally
geodesic Seifert surface, which is totally knotted, since the complement
of a regular neighborhood
of the  surface is Thurston's wye manifold, which is known to be hyperbolic
with incompressible boundary.

\medskip

\section{Knots without Totally Geodesic Seifert
Surfaces}\label{S:twobridge}

We prove the following theorem.

\begin{thm}
Hyperbolic two-bridge knots do not possess orientable totally geodesic
surfaces.
\end{thm}

\begin{proof}
Note that the only non-hyperbolic $2$-bridge knots are the $(2,q)$-torus
knots. There are no non-boundary parallel incompressible closed
surfaces in their complements.  In ~\cite{HT85}, a complete 
classification of all the
incompressible boundary-incompressible surfaces in the complements of
two-bridge knots is given.  All such surfaces are carried by branched
surfaces that appear as in Figure \ref{fig:branchedtwobridge}. Such a
surface is
orientable if and only if each $b_i$ is even
where $p/q = r +[b_1,b_2, \ldots, b_n]$. For each $p/q$, there is a unique
choice of such even $b_i$'s.

\begin{figure}[htbp]
\begin{center}
\includegraphics*[height=3.5 in]{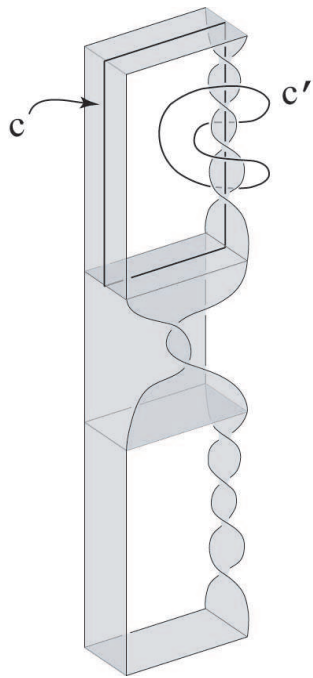}
\end{center}
\caption{\label{fig:branchedtwobridge}All Seifert surfaces for 2-bridge
knots take this basic form.}
\end{figure}

Note that a Seifert surface corresponds to  a $1$-sheeted surface 
carried by the
branched surface. Each possible Seifert surface is
obtained as a collection of horizontal "plumbing" disks either to the
inside or the outside and vertical strips between them, each with an even
number of half-twists. The other incompressible boundary-incompressible
  surfaces have more than one sheet.

If for some $i= 1, 2, \dots, n$, we have $b_i=2$, then there is a disk in
$\Sth$ with boundary that consists of two arcs on $S$ and two arcs on the
cusp
boundary. When the
manifold is doubled along $S$, this disk doubles
to an essential annulus in the resulting manifold.
Hence, the resulting manifold cannot be hyperbolic, implying that the
surface $S$ could not be totally geodesic.

When $b_i > 2$ for all i, choose a curve $c$ on $S$ as in Figure
\ref{fig:branchedtwobridge}. Note that the curve
$c$ can be isotoped to $c'$ as in Figure \ref{fig:branchedtwobridge}. The
curve $c'$ sits
on a torus $T$ in $\Sth - N(S) - N(K)$. Cutting the torus open along $c'$
yields an annulus in $\Sth - N(S) - N(K)$, both
boundary components of which lie on $\partial N(c')$ in $T$. Gluing on to
its boundary components a pair of parallel annuli that give the isotopy
from $c$ to $c'$ yields
an annulus $A$ in $\Sth-N(S)$ with boundary in $\partial N(S)$.  The
annulus $A$ is incompressible in $\Sth-int(N(S))$ since its core curve is
a power  of a generator for the fundamental group of the
complementary handlebody. To see that it is boundary-incompressible, we can
choose the essential arc $\alpha$ in $A$ shown in Figure
\ref{fig:essentialarc}.  If
$\alpha$ can be isotoped into the surface, then there exists a disk with
boundary consisting of two arcs on $S$ and two arcs on the cusp boundary.
Again, when the manifold is doubled along $S$, this yields an essential
annulus, a contradiction to the hyperbolicity of a hyperbolic manifold
doubled along a totally geodesic surface within it.  If $\alpha$
cannot be isotoped into the surface, then $A$ is an essential annulus
intersecting $S$ in its boundary, which obstructs $S$ from being totally
geodesic.

\begin{figure}[htbp]
\begin{center}
\includegraphics*[height=2 in]{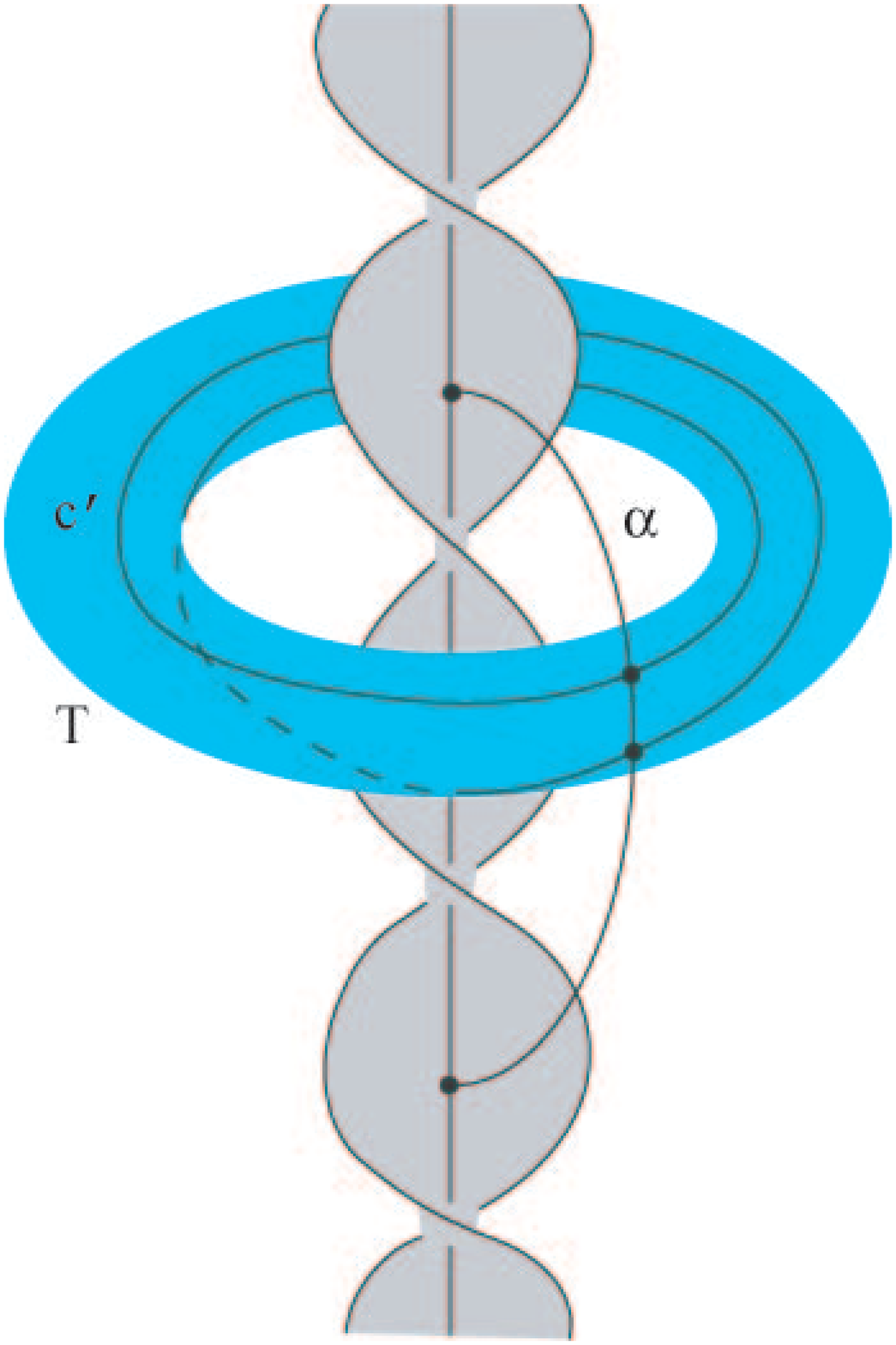}
\end{center}
\caption{\label{fig:essentialarc}}
\end{figure}

\end{proof}

Do two-bridge knots possess non-orientable totally geodesic surfaces?
Any essential non-orientable surface for a
2-bridge knot or link is given by  a choice of $[b_1,b_2, \ldots, b_n]$,
where $p/q = r +[b_1,b_2, \ldots, b_n]$
and at least one $b_i$ is odd. The only
candidates for totally geodesic
   surfaces from amongst these are those
where each $b_i$ is odd. In the case
   of links, these surfaces can be
totally geodesic. The Whitehead link
   possesses such a surface. However,
we have yet to see such a surface in a
   2-bridge knot complement that is
totally geodesic.

\end{document}